\documentclass[12pt,a4paper]{amsart}

\usepackage{amscd,amssymb,amsthm}

\title[constants in Jackson--Stechkin Theorem]
{Bohr--Favard Inequality for differences and
constants in Jackson--Stechkin Theorem}

\author{Y. Kryakin}
\address{Institute of Mathematics\\
University of Wroc{\l}aw \\
pl. Grunwaldzki 2/4 \\
50-384 Wroc{\l}aw, Poland}
\email{kryakin@math.uni.wroc.pl}

\begin{document}
\maketitle

\begin{abstract}
We consider uniform approximations by trigonometric
polynomials. The aim  of the paper is to obtain good estimates of
the Jackson--Stechkin constants $J_m$. We prove that $ J_m \le C
2^{-m+5/2\log_2m}$. Our proof is based on the difference analogue
of the Bohr--Favard inequality.
\end{abstract}

\vskip 1cm

\section{Introduction}\label{sec1}

Let $C(T)$ be the space of continuous function on the  unit
circle $T=[-\pi,\pi)/2\pi Z$.  Denote by $ T_n $ the space of
trigonometric polynomials $ \sum_{j=-n}^{n} c_j e^{ijt}$  of
degree $ \le n$. Let  $ C(T_n^\perp )$ be the class of functions
from $C(T)$ which have no spectrum in $[-n,n]$.

We are interested in the uniform approximations of $f \in { C(T)}$
by trigonometric polynomials. Set
$$
 \| f \| := \max_{x \in T} |f(x)|, \quad \|f \|_1 := \int_T
|f(u)| \, du.
$$
We shall denote by
$$
E_n(f):= \inf_{\tau \in
T_n} \| f - \tau \|
$$
the best approximation of $ f \in C(T)$ by $ \tau \in T_n$.
\noindent Define  $m$--difference of $f$ at a point $x$ with
step  $h$ as follows
$$ \Delta_h^m f(x) :=
\sum_{i=0}^m (-1)^{m-i} \binom m i f(x-(m/2)h +ih).
 $$
The $m$--th modulus of smoothness of $f$ is defined by $$ \omega_m
(f, \delta):=\sup_{x, \ |h|<\delta} |\Delta_h^m f(x)|. $$ The
following theorem is one of the central results in approximation
theory.
 \vskip .2cm
 { \bf Theorem J.} {\it If $f \in
C(T)$, then

$$ E_n(f) \le J_m \omega_m (f, \frac{2\pi}n), \quad n \ge m. $$ }

Theorem J was proved by  Jackson \cite{J} ($m=1$),
Akhiezer  \cite{A} ($m=2$), Stechkin \cite{S} \ ($m>2$).
Stechkin's proof \cite[p.205]{DL}
 does not allow us to obtain the estimate
 $J_m <  m^m$ for  $m>m_0$. We  shall
 show that constant $J_m$ tends to zero exponentially
 as $m \to \infty$.
The main result of the paper is the following.

\vskip .2cm
{ \bf Theorem 1.} {\it For $f \in C(T), \  \alpha := 2^{-5/2} 3^{1/2} \pi <
1$,

$$ E_n(f) \le
\frac{2(k+1)^2}{\binom{2k}k} \ \omega_{2k} (f,
\frac{2\pi}{ n} \alpha (1+1/k) ), \quad n \ge 2k. \eqno (1) $$ }

\vskip .2cm
{ \bf Corollary.} { $$ J_m \le C \ 2^{- m  + (5/2)
\log_2 m}. \eqno (2) $$ }

\vskip .1cm \noindent There is a connection between (1) and the
following theorem  for functions with uniform oscillation on
$[a,b]$. We write  $f \in O_n (a,b)$ if $f \in C(a,b)$ and

$$ \int_{a}^{a+i(b-a)/n} f(u) \, du =0, \quad  i=1,\dots,
n.
$$

\vskip .2cm

{ \bf Theorem  W.} {\it For $f \in  O_n(a,b), $
$$
\| f \| \le W_n \omega_m (f, \frac {b-a}n), \quad n \ge m. $$ }

\noindent
It is known that $W_n <2 + 1/e^2$ and $W_n =1$ for $n<8$
\cite{GKS,Z}. Note that our conjecture  $W_n=1$ implies Sendov's
 conjecture  \cite{S1}.

\vskip .2cm

\noindent
The condition $f \in O_n (a,b)$  is similar to the condition $f \in
C(T_n^\perp)$. The Bohr--Favard \cite{B,F} inequality reads

\vskip .2cm { \bf Theorem F.} {\it If  $f, \ f^{(m)}  \in  C (T_n^\perp)$,
then  $$ \| f \| \le F_m (n+1)^{-m} \| f^{(m)} \|, \qquad F_m: =
\frac{4}\pi\sum_{i=0}^\infty (-1)^{i(m+1)} (2i+1)^{-m-1}. $$ }
\noindent
The following difference analogue of Theorem F
 is the key result of the paper.

\vskip .2cm

{\bf Theorem 2.} {\it If  $f \in C (T_n^\perp )$, then

$$ \| f \| \le \frac{k+1}{\binom{2k}k} \omega_{2k}
(f,  \frac{2\pi}{n} \alpha), \quad n \ge 2k. $$ }

\noindent
The third theorem in  this paper gives a link between Theorem 2
and Theorem 1. It is devoted to the approximation by
Vall\'{e}e Poussin means.
\noindent
Denote by $s_i $ the operator of  $i$--th partial Fourier sum.
Let
$$ v_{k,m}:= \frac 1m \sum_{i=km}^{(k+1)m-1} s_i. $$

\vskip .2cm

{\bf Theorem 3.} {\it If $f \in C (T)$, then
$$ \| f-v_{k,m}f \|
\le  \frac{2(k+1)^2}{\binom{2k}k }
 \omega_{2k} (f,  \frac{2\pi}{km}\alpha), \quad m \ge 2 .
$$ }

\noindent
We shall use  special pointwise difference operators $w^s_{x,h,2k}$.
The operators $w_{x,h,2k}^1$ were introduced by Ivanov and Takev \cite{IT}
under the influence of Beurling's proof of Whitney's theorem in $R$ and
$R_+$ \cite[p.83]{W}.
We should also mention  Brudnyi's paper \cite{B1}. It contains the
construction of  smoothing operators,
which is widely used in approximation theory
\cite[p.177]{DL}. The operators
$w^s_{x,h,2k}$ are similar to Brudnyi's operators, but provide
more delicate estimates.

\section{Smoothing operators}\label{sec2}

Let $I_x$  be the identity  operator at the point $x \in T$.
$$
I_x f:= f*\delta_x = f(x).
$$
\noindent
Averaging operator on  $[x-h, x+h ], \  0<h<\pi, $  will be denoted by  $I_{x,h}$.

$$ I_{x,h}f:= \frac{1}{2h} \int_{x-h}^{x+h} f(u) \, du. $$
\noindent
It is clear that
$$
I_{x,h}f = \int_T f(u) B_h^1 (x-u) \, du = (f*B_h^1)(x),
$$
where
$$ B_h^1 (x):= \begin{cases} \frac 1{2h}, \quad & x \in
[-h, h]\\ 0 \quad &x \in [-h, h]^c.
\end{cases}
$$
\vskip .2cm
\noindent
Set $I_{x,h}^0:=I_x, \ I_{x,h}^1:=I_{x,h}.
$ Define $I_h^s, \ s=2, \dots , $ by

$$
I_{x,h}^s f :=(I_{x,h})^s f = (f*B_h^s)(x),
\quad B_h^s (x):=(B_h^1*B_h^{s-1})(x).
$$
\vskip .1cm
\noindent
Note that the nonnegative  function  $B_h^{s}$  has the following
properties:

$$ \int_T B_h^s (u) \, du =1, $$

$$ \mbox{  supp\, } B_h^s (u) = [-sh, sh], \quad 0<sh<\pi. $$
\vskip .2cm
\noindent
Operator of differentiation $D$ acts on $I_{x,h}$ in the following way:

$$ DI_{x,h}^s =(2h)^{-1}(I_{x+h,h}^{s-1} - I_{x-h,h}^{s-1}). $$
\noindent
This implies

$$
D^m I_{x,h}^s = (2h)^{-m} \sum_{i=0}^m (-1)^{m-i} \binom mi I^s_{x-mh+2ih}
, \quad s \ge m.
$$
\vskip .2cm
\noindent
The operators $w_{x,h,2k}^s$ measuring the local
$2k$--th smoothness of a function are defined by

$$ w_{x,h,2k}^s := \sum_{i=0}^{2k} (-1)^{k-i}
\frac{\binom{2k}{i}}{\binom{2k}k } I_{x,(i-k)h}^s =I_x - 2
\sum_{i=1}^{k} (-1)^{i+1}
 \frac{\binom{2k}{k+i}}{\binom{2k}k } I_{x,ih}^s,
 \quad I_{x,-ih}^s:=I_{x,ih}^s. \eqno (3)
$$
\noindent
We can write $w_{x, h, 2k}^s f$ as

$$
w_{x,h,2k}^s f =(-1)^k \binom{2k}k^{-1} \int_T \Delta_u^{2k} f(x)B_h^s(
u) \, du, \quad 0<skh<\pi. $$
\noindent
Note that

$$ \| w_{*,h,2k}^s f\|
\le \binom{2k}k^{-1} \omega_{2k}(f,sh).
$$

\section{Proof of Theorem 2}\label{sec3}

Rewrite (3) in the following form

$$ I_x = w_{x,h,2k}^s +  2 \sum_{i=1}^{k} (-1)^{i+1}
 \frac{\binom{2k}{k+i}}{\binom{2k}k } I_{x,ih}^s. \eqno (4)
$$

\noindent We shall show that for  $f \in  C(T_n^\perp) $ one can
choose $s,h$ \ such that

$$ \| I_{*,ih}^s  f \| \le c(h,s,i) \| f \|, $$
and
$$ 2
\sum_{i=1}^{k} \frac{\binom{2k}{k+i}}{\binom{2k}k } c(h,s,i) \le
b(k) < 1. $$

\noindent
From
$$ I_{x,ih}^s f = (B_{ih}^s * f) (x), $$
we have for    $f \in
C(T_n^\perp),  \ \tau \in T_n $,

$$ I_{x,ih}^s f = ((B_{ih}^s - \tau )*f)(x). $$
This gives
$$ \| I_{*,ih}^s f\| \le \inf_{\tau \in
T_n} \|B_{ih}^s -
\tau  \|_1 \|f\|.
$$
\vskip .1cm
\noindent
By  Favard--Nikolskii theorem \cite[p.215]{DL}
$$ \inf_{\tau \in T_n} \| B_{ih}^s -
\tau \|_1
 \le F_{s-1} n^{-s+1} \| D^{s-1} B_{ih}^{s} \|_1.
$$
The equality
$$ D^{m} B_{h}^s(x) = (2h)^{-m} \Delta_{2h}^s B_h^{s-m}(x),
\quad s \ge m, $$
implies
$$\inf_{\tau \in T_n}  \|B_{ih}^s - \tau \|_1 \le
F_{s-1} (2ihn)^{-s+1} \| \Delta_{2ih}^{s-1} B_{ih}^1
\|_1. $$
Choose  $s=3, \  \beta:=\frac{\pi}{\sqrt{6}},  \ h= \beta \pi
(2n)^{-1}  $. We thus get
$$ c(h,3,i) \le \inf_{\tau \in T_n} \|B_{ih}^s - \tau  \|_1
\le 4 F_2 {( \beta \pi i)}^{-2} = \frac 1{2 \beta^2 i^2}, $$
and
$$ 2 \sum_{i=1}^{k}
\frac{\binom{2k}{k+i}}{\binom{2k}k } c(h,3,i) \le
\frac 1{\beta^2} \frac{\binom{2k}{k+1}}{\binom{2k}k }\left(
\sum_{i=1}^{k} i^{-2} \right)   \le
\frac{\binom{2k}{k+1}}{\binom{2k}k }<1.
$$
\vskip .1cm
\noindent
The identity (4) at the extremal point $x_0$, such that \ $|f(x_0)| = \|f \|,$ gives

$$ \| f \| (1-\frac{\binom{2k}{k+1}}{\binom{2k}k }) \le \|
w_{*,h,2k}^3f)\|. $$
The inequality
$$ \| w_{*,h,2k}^3 f \| \le
\frac{1}{\binom{2k}k} \omega_{2k} (f,3h), $$
implies
$$ \| f \| \le \frac{1}{\binom{2k}{k} -\binom{2k}{k+1}}
\omega_{2k}(f, \frac{3 \beta \pi}{2n})= \frac{k+1}{\binom{2k}{k}}
\omega_{2k}(f,  \frac{2\pi}{n} \alpha). $$
\noindent
This completes the proof of Theorem 2.

\section{Proofs of Theorem 3 and Theorem 1}\label{sec4}

Note that   Vall\'{e}e Poussin means $v_{k,m}$
are the simple combination of Fejer's means.

$$ v_{k,m}= \frac 1m \sum_{i=km}^{(k+1)m-1} s_i
=(k+1)\sigma_{(k+1)m-1} - k\sigma_{km-1},
$$
where
$$
\sigma_j:=\frac{1}{j+1} \sum_{i=0}^{j} s_i. $$
For Fejer's means $\sigma_j$ we have
$$ \|\sigma_jf \| \le \|f \|,
\quad \Delta_h^{2k} \sigma_jf = \sigma_j \Delta_h^{2k} f. $$
Therefore
$$
 \omega_{2k} (v_{k,m}f ,h)  \le (2k+1) \omega_{2k} (f ,h).
$$
\vskip .1cm
\noindent
Let $g:=f-v_{k,m}f \in C(T_{km}^\perp)$.  Theorem 2
implies
$$ \| g \| \le \frac{k+1}{\binom{2k}{k}}
\omega_{2k} (g,  \frac{2\pi}{km} \alpha) \le
\frac{2(k+1)^2}{\binom{2k}{k}
}  \omega_{2k} (f,  \frac{2\pi}{km} \alpha). $$
Theorem 3 is proved.
\vskip .2cm
\noindent
For the best approximation by trigonometric polynomials
of degree $n =km+i, \ i \in [0, m-1], \ m \ge 2,$ we get
$$ E_n(f)
\le E_{km}(f) \le \| f-v_{m,k}f \| \le \frac{2(k+1)^2}{\binom{2k}{k}
}  \omega_{2k} (f, \frac{2\pi}{km} \, \alpha). $$
From
$$ n/(km) \le (km+m-1)/(km) <
1+1/k, $$
we have
$$ E_n(f) \le \frac {2(k+1)^2}{\binom{2k}{k} }
\omega_{2k} (f, \frac{2\pi}{n} \, \alpha (1+1/k)).
$$
This proves Theorem 1.
\vskip .2cm
\noindent
The inequalities  $$ \omega_{2k} (f, \delta)  \le
2 \omega_{2k-1}(f,\delta), \qquad \frac{4^k}{\sqrt{\pi(k+1/2)}}<
\binom{2k}k < \frac{4^k}{\sqrt{\pi k}}, $$
lead to (2).

\end{document}